\documentclass[12pt]{article}

\usepackage{epsfig}
\usepackage{amssymb}
\usepackage{amsmath}

\def\shadowbox{\hbox{\rule[-0.0ex]{0.1ex}{1.2ex}%
\hspace{-0.1ex}\rule[-0.0ex]{1.2ex}{0.1ex}%
\hspace{0.0ex}\rule[-0.0ex]{0.1ex}{1.2ex}\hspace{-1.3ex}%
\rule[1.15ex]{1.25ex}{0.1ex}\hspace{-0.0ex}\rule[-0.25ex]{0.3ex}{1.1ex}%
\hspace{-1.2ex}\rule[-0.25ex]{1.1ex}{0.25ex}}}
\def\qed{\ifmmode \hbox{\hfill\shadowbox}
     \else \hphantom{x}\hfill\shadowbox \fi}

\newtheorem{theorem}{Theorem}[section]
\newtheorem{lemma}[theorem]{Lemma}
\newtheorem{definition}[theorem]{Definition}

\newtheorem{corollary}[theorem]{Corollary}

\def\remark{{\noindent \bf Remark:\hspace{0.5em}}}

\def\eps{\varepsilon}
\def\Cst{{\mathbb C}}
\def\Nst{{\mathbb N}}
\def\Qst{{\mathbb Q}}
\def\Rst{{\mathbb R}}

\def\Zst{{\mathbb Z}}
\def\Bsp{{\boldsymbol B}}
\def\Lsp{{\boldsymbol L}}
\def\lsp{{\boldsymbol\ell}}
\def\LtR{{\Lsp_2(\Rst)}}
\def\ltZ{{\lsp_2(\Zst)}}
\def\range{\operatorname{range}}
\def\Image{\operatorname{Im}}
\def\cond{\mbox{cond\/}}
\def\ltZt{{\lsp^2(\Zst^{2})}}
\def\ninN{{n{\in}\Nst}}
\def\compact{{\cal C}_{[\alpha,\beta]}}
\def\compacti{{\cal C}_{[\alpha_1,\beta_1]}}
\def\expmat{{\cal E}_{\lambda}}
\def\expmati{{\cal E}_{\lambda_1}}
\def\polmat{{\cal Q}_{s}}

\newcommand{\gkl}{g_{k/b,l/a}}
\newcommand{\psikl}{\psi_{k,l}}
\newcommand{\gklp}{g_{k'/b,l'/a}}

\newcommand{\gamman}{\gamma^{(n)}}
\newcommand{\gtight}{\tilde{g}}
\newcommand{\gtightkl}{\gtight_{k/b,l/a}}
\newcommand{\gtightnm}{\gtight_{na,mb}}
\newcommand{\stft}{{\cal V}}
\newcommand{\gnm}{g_{na,mb}}

\newcommand{\gnmp}{g_{n'a,m'b}}
\newcommand{\gab}{(g,a,b)}
\newcommand{\gba}{(g,1/b,1/a)}
\newcommand{\ckl}{c_{k,l}}

\newcommand{\cnm}{c_{n,m}}
\newcommand{\dnm}{d_{n,m}}
\newcommand{\Pn}{P_n}
\newcommand{\Hn}{H_n}
\newcommand{\Ha}{H^{\ast}}
\newcommand{\HH}{H\Ha}
\newcommand{\Hna}{\Ha_{n}}
\newcommand{\HHn}{\Hn \Hna}
\newcommand{\HHI}{(H\Ha)^{-1}}
\newcommand{\HHnI}{(\HHn)^{-1}}
\newcommand{\Tn}{T_n}
\newcommand{\Tna}{T_n^{\ast}}

\newcommand{\cn}{c^{(n)}}
\newcommand{\dn}{d^{(n)}}
\newcommand{\xn}{x^{(n)}}

\newcommand{\wienerw}{{\cal A}_w}
\newcommand{\laurentw}{{\cal W}_w}
\newcommand{\wienerwm}{{\cal A}_w^{m}}
\newcommand{\laurentwm}{{\cal W}_w^{m}}
\newcommand{\sigman}{\sigma^{(n)}}
\def\ntoinf{{n \rightarrow \infty }}
\def\toinf{{\rightarrow \infty }}
\def\tozero{{\rightarrow 0 }}
\newcommand{\gamnm}{\gamma_{n,m}}

\begin{document}

\title{Approximation of dual Gabor frames, window decay, and wireless
communications}

\author{Thomas Strohmer\thanks{Department of Mathematics, University of 
California, Davis, CA 95616-8633, USA; Email: strohmer@math.ucdavis.edu. 
This work was supported by NSF grant 9973373.}}

\date{}
\maketitle
\vspace*{-6cm}
\noindent
{\small \tt Submitted to Applied and Computational Harmonic Analysis}

\vspace*{6cm}

\begin{abstract}
We consider three problems for Gabor frames that have recently
received much attention. The first problem concerns the approximation
of dual Gabor frames in $\LtR$ by finite-dimensional methods. Utilizing
Wexler-Raz type duality relations we derive a method to approximate 
the dual Gabor frame, that is much simpler than previously proposed
techniques. Furthermore it enables us to give estimates for the approximation
rate when the dimension of the finite model approaches infinity.
The second problem concerns the relation between the decay of the 
window function $g$ and its dual $\gamma$. Based on results on
commutative Banach algebras and Laurent operators we derive a general
condition under which the dual $\gamma$ inherits the decay properties
of $g$. 
The third problem concerns the design of pulse shapes for orthogonal frequency 
division multiplex (OFDM) systems for time- and frequency dispersive
channels. In particular, we provide a theoretical foundation for a
recently proposed algorithm to construct orthogonal transmission
functions that are well localized in the time-frequency plane.
\end{abstract}



\section{Introduction}
\label{s:intro}

Gabor systems play an important role in signal processing and
digital communication. In filter bank theory they are known under the
name {\em oversampled modulated filter banks}~\cite{CR83}, in wireline
communications they correspond to the concept of 
{\em discrete multitone transmultiplexing}, and in wireless
communications they are (implicitly) used in
{\em orthogonal frequency division multiple access
systems}~\cite{Zou95,KM98,Bol99}.

A Gabor system consists of functions of the form
\begin{equation}
\gnm(t) = e^{2\pi i mb} g(t-na),\,\,\, n,m \in \Zst,\,\, a,b \in \Rst,
\label{gaborsystem}
\end{equation}
where $g \in \LtR$ is -- depending on the context -- called {\em window}, 
{\em atom}, or {\em pulse shape}.
The parameters $a$ and $b$ represent the time-shift and frequency-shift,
respectively.

We say that $(g,a,b)$ generates a {\em Gabor frame} for $\LtR$ for given shift 
parameters $a,b$ if there exist constants ({\em frame bounds}) $A,B>0$ such 
that
\begin{equation}
A \|f\|^2 \le \sum_{n,m \in \Zst} |\langle f, \gnm \rangle|^2 \le B\|f\|^2,
\label{gaborframe}
\end{equation}
for any $f \in \LtR$.

The analysis operator $T$ is defined as
\begin{equation}
\label{analysis}
T : f \in \LtR \rightarrow Tf = \{\langle f, \gnm \rangle\}_{n,m \in \Zst},
\end{equation}
and the synthesis operator, which happens to be the adjoint of $T$ is
\begin{equation}
\label{synthesis}
T^{\ast} : c \in \lsp_2(\Zst \times \Zst) \rightarrow 
T^{\ast} c = \sum_{n,m\in \Zst} c_{nm} \gnm .
\end{equation}

The Gabor frame operator is defined by
\begin{equation}
Sf = \sum_{n,m \in \Zst} \langle f, \gnm \rangle \gnm , \qquad f \in \LtR,
\label{frameop}
\end{equation}
and satisfies 
$$ I A \le S \le I B$$
where $I$ is the identity operator on $\LtR$. Of course $S = T^{\ast} T$.

If $(g,a,b)$ establishes a Gabor frame for $\LtR$ then any $f$ in $\LtR$ 
can be represented as
\begin{equation}
\label{framerep}
f = \sum_{n,m \in \Zst} \langle f,\gamnm \rangle \gnm = 
\sum_{n,m \in \Zst} \langle f,\gnm \rangle \gamnm, 
\end{equation}
where the {\em dual frame} $\{\gamnm\}$ is given by
$$\gamnm = e^{2\pi i mb} \gamma(t-na),\,\,\, n,m \in \Zst;\, a,b \in \Rst,
$$
with $\gamma = S^{-1} g$. 
In general there are many functions generating dual frames that satisfy 
relation~\eqref{framerep}. The ``canonical'' dual window $\gamma$ has several 
nice properties. One of them is that it has minimal $\Lsp_2$-norm among
all dual functions. In this paper we concentrate on the canonical
dual window and henceforth simply talk about {\em the} dual window
and {\em the} dual frame. For more details about the properties
of Gabor frames and their duals the reader is referred to~\cite{Dau90,FS98}.

The rest of the paper is organized as follows. 
In Section~\ref{s:finite} we analyze the problem of approximating the
dual window by using finite-dimensional methods. Based on 
Wexler-Raz type duality relations we derive a method to approximate 
the dual Gabor frame, that is much simpler than previously proposed
techniques. Furthermore we show that for windows with exponential decay
in time and frequency the proposed approach yields an exponential
approximation rate when the dimension of the finite model approaches infinity.
In Section~\ref{s:laurent} we dig deeper into decay properties of
$g$ and its dual $\gamma$. Based on results on commutative Banach algebras 
and Laurent operators we derive a general condition on the decay of $g$
which guarantees that the dual $\gamma$ inherits these decay properties.
Finally in Section~\ref{s:ofdm} we demonstrate the relevance 
of the results derived in Section~\ref{s:laurent} for wireless communications.
In particular, we provide a theoretical foundation for a
recently proposed algorithm to construct an orthogonal frequency
division multiplex (OFDM) system with good time-frequency localization
properties.

\bigskip
Before we proceed we introduce a few notations used throughout the paper.

The Fourier transform of a function $f$ is given by
$$
\hat{f}(\omega)=\int\limits_{-\infty}^{+\infty} f(t) e^{-2\pi it\omega} dt.
$$
The {\em short time Fourier transform} (STFT) of $f$ with respect to the 
(sufficiently nice) window $g$ is
\begin{equation}
(\stft_g f)(t,\omega) = \int \limits_{-\infty}^{+\infty} 
f(x) \overline{g(x-t)} e^{-2\pi i x\omega} \, dx. 
\label{stftdef}
\end{equation}

A locally integrable function $w$ is called {\em weight function}, if $w$ 
is positive and submultiplicative, i.e., if $w(t)>0$ and 
$w(t_1+t_2) \le w(t_1) w(t_2)$.
 
The space $\Lsp_{1,w}(\Rst)$ consists of all functions $f$ with
$$\int \limits_{-\infty}^{+\infty} |f(t)| w(t) dt < \infty ,$$
where $w$ is a weight function. Similarly $\lsp_{1,w}(\Zst)$
consists of all sequences $x=\{x_k\}_{k \in \Zst}$ with
$$\sum_{k = \infty}^{\infty}|x_k| w(k) < \infty.$$

It is convenient to define following spaces:
$$\compact=\Big\{f \in \LtR :\,\exists\, \alpha, \beta \,\,\text{such that}\,\,
|f(t)|=0 \,\,, \forall |t|\notin [\alpha, \beta]\Big\}.$$
$$\expmat =\Big\{ f \in \LtR :\, \exists\, \lambda>0, c>0 \,\,
\text{such that}\,\, |f(t)|\le c e^{-\lambda |t|}\,\,, \forall |t|\Big\}.$$
$$\polmat =\Big\{ f \in \LtR :\, \exists\, s>1, c>0 \,\,
\text{such that}\,\, |f(t)|\le c (1+|t|)^{-s}\,\,, \forall |t|\Big\}.$$

Finally, the Moore-Penrose inverse~\cite{EHN96} of a bounded operator 
$T$ is denoted by $T^{+}$.

\section{Gabor frames, finite sections and the duality condition} 
\label{s:finite}

The theoretical concepts for Gabor analysis are usually developed
for infinite-dimensional function spaces, notably for $\LtR$ or $\ltZ$, 
whereas all numerical implementations have to be done within a 
finite-dimensional framework. The connection between
Gabor systems on $\ltZ$ and several finite-dimensional models has been
clarified in~\cite{Str98a,Str99}. In this section we extend these
results to Gabor frames for $\LtR$. Different techniques than those
used for $\ltZ$ are required for $\LtR$.

Our approach relies on a remarkable property of Gabor frames,
whose discovery has its origin in a paper by Wexler and Raz~\cite{WR90},
their result was later made precise and extended by Janssen~\cite{Jan95},
Ron and Shen~\cite{RS97a}, and Daubechies, H.~Landau, and 
Z.~Landau~\cite{DLL95}.

For given $g, a, b$ we define the operator $H$ by 
$$ Hf :=\{ \langle f , \gkl  \rangle \}_{k,l \in \Zst}, \qquad f \in \LtR,$$
where
$$\gkl(t) = g(t-k/b)e^{2\pi i tl/a}.$$
The adjoint $\Ha$ of $H$ is 
$$\Ha c = \sum_{k,l} \ckl \gkl , \qquad c=\{\ckl\}_{k,l \in \Zst} \in 
\ltZt.$$

We identify $H \Ha$ with its matrix representation, given by
$$H \Ha = \{\langle \gklp , \gkl \rangle \}_{k,l;k',l' \in \Zst}. $$

The Wexler-Raz biorthogonality relations imply that the dual window
satisfies~\cite{Jan95,DLL95}:
\begin{equation}
\label{wexlerraz}
\gamma= \Ha (H \Ha)^{-1} \sigma = ab \sum_{k,l \in \Zst} 
[(H\Ha)^{-1}]_{k,l;0,0} \,\, \gkl,
\end{equation}
where $\sigma = \{(ab) \delta_{k0} \delta_{l0} \}_{k,l \in \Zst}$.

Furthermore, there holds:
\begin{theorem} [\cite{RS97a}]
\label{ronshen1}
$\gba$ is a Riesz basis for its closed linear span if and only if $\gab$
is a frame for $\LtR$.
\end{theorem}

Formula~\eqref{wexlerraz} and Theorem~\ref{ronshen1} are the main 
ingredients for our approach to approximate the dual window $\gamma$ using a 
finite-dimensional model.

Now, for $x \in \ltZt$ and $n \in \Nst$ define the orthogonal projections
$\Pn$ by
$$(\Pn x)_{k,l} = 
\begin{cases}
x_{k,l} & \text{if $\max\{|k|,|l|\} \le n$}, \\
0       & \text{else.}
\end{cases}
$$
We identify the image of $\Pn$ with the $(2n+1)^2$-dimensional space
$\Cst^{(2n+1)^2}$ and write
$$\Hn f := \Pn Hf = \{\langle f , \gkl \rangle\}_{|k|,|l| \le n}.$$
The matrix 
$$\Hn \Hna = (H \Ha)_n = \Pn H \Ha \Pn = 
\{\langle \gklp , \gkl \rangle \}_{|k|,|l|;|k'|,|l'| \le n}. $$
is a finite section of the infinite-dimensional Gram matrix $H \Ha$.

We say that the finite section method is {\em applicable} to $H \Ha$,
if, beginning with some $n \in \Nst$, for each $y \in \range (H)$ the equation
$$\Hna [\Hn \Hna]^{-1} \xn = \Pn y$$
has a unique solution $\xn \in \Image (\Pn)$ and as $\ntoinf$ the vectors
$\xn$ tend to the solution of $H \Ha x = y$.

We set
$$\gamman := \Hna [\Hn\Hna]^{-1} \sigman$$
for $n=0,1,2,\dots$. 

There holds:
\begin{theorem}
\label{th1}
Let $\gab$ generate a Gabor frame for $\LtR$ and let $\gamma$ be the dual
window. Then
$$\|\gamma - \gamman\| \tozero \qquad \text{for $\ntoinf$}.$$
\end{theorem}

We need following result for the proof of Theorem~\ref{th1}.

\begin{lemma}[Lemma of Kantorovich, \cite{RM94}]
\label{kantorovich}
Let $\{K_n\}_{n\in \Nst}$ be a sequence of invertible operators on a Banach
space $\Bsp$, and assume that the sequence is uniformly bounded above
and below, i.e, that there exist constants $0<C_1, C_2<\infty$ such that
\begin{equation}\label{eq:unifbound}
C_1 \, \|f\|_\Bsp \leq \|K_n f\|_\Bsp \leq C_2 \, \|f\|_\Bsp
  \qquad\text{for all $\ninN$ and $f\in\Bsp$.}
\end{equation}
Then $\{K_n\}$ converges in the strong operator topology
to an invertible operator if and only if the same is true for
$\{K_n^{-1}\}_n$, and then
$$
  \{\lim {K_n}\}^{-1} = \lim \{K_n^{-1}\} \,.
$$
\end{lemma}

\noindent
{\bf Proof of Theorem \ref{th1}:}

Theorem~\ref{ronshen1} implies that $\gba$ is a Riesz basis for
its closed span with {\em Riesz bounds} $A, B$. Any finite subset of a 
Riesz basis is again a Riesz basis for its closed span, cf.~Chapter 1 
in~\cite{You80}.  Thus for $\{\gkl\}_{|k|,|l| \le n}, n \in \Nst$ there 
exist constants $A_n, B_n$ such that 
\begin{equation}
A \le A_n \le \Hn\Hna \le B_n \le B.
\label{rieszbounds}
\end{equation}
Hence we can apply Lemma~\ref{kantorovich} to establish that the
finite section method is applicable to $H \Ha$.

Finally note that $\Hn$ obviously converges to $H$ pointwise for $\ntoinf$,
thus $\|\gamman - \gamma\| \rightarrow 0$ for $\ntoinf$.  \qed

It is interesting to ask if we can give some estimate on the
rate of approximation of the finite section method in case the window
$g$ satisfies certain decay conditions in time and/or frequency.
In the following theorem we concentrate on windows with exponential decay 
in time and frequency domain, but it will become clear from the proof that 
similar results may be obtained for other types of decay.

\begin{theorem}
Let $\gab$ be a Gabor frame for $\LtR$ and assume that there exist
constants $C,D>0$ such that
\begin{equation}
\label{expdecay}
|g(t)| \le C e^{-\lambda |t|}, \quad \text{and} \quad 
|\hat{g}(\omega)| \le D e^{-\lambda |\omega|},
\end{equation}
then there exists a $\lambda' < \lambda$ and a constant $C'$ depending 
on the frame bounds and on $\lambda'$, but independent of $n$ such that
$$ \|\gamma - \gamman \| \le C' e^{-\lambda' n}. $$
\end{theorem}

\begin{proof}
First we show that~\eqref{expdecay} implies that
the entries of $H\Ha$ satisfy
\begin{equation}
\label{matexp}
(H\Ha)_{kl;k'l'} \le D_2 e^{-\lambda_1 (|k-k'|+|l-l'|)}
\end{equation}
for some $\lambda_1 < \lambda$ and some constant $D_2$.

It is clear that 
\begin{equation}
\label{innprod}
|(H\Ha)_{kl,k'l'}|=|\langle \gkl, \gklp \rangle | = 
|\langle g, g_{k-k',l-l'} \rangle |
\end{equation}
and
\begin{equation}
\label{sampledstft}
\langle g, \gkl \rangle = (\stft_g g)(k/b,l/a) \quad k,l \in \Zst.
\end{equation}

Now
\begin{gather}
|(\stft_g g)(t,\omega)| = \Bigl| \int \limits_{-\infty}^{+\infty}  
g(x) \overline{g(x-t)} e^{-2\pi i \omega x}\Bigr| dx \le
\int \limits_{-\infty}^{+\infty} |g(x)| |\overline{g(x-t)}| dx \notag \\
\le C^2\int\limits_{-\infty}^{+\infty} e^{-\lambda |x|}e^{-\lambda |x-t|}dx 
\le C^2\int\limits_{-\infty}^{+\infty} e^{-\lambda |x|}e^{-\lambda_1 |x-t|}dx 
\end{gather}
for some $\lambda_1 < \lambda$. Set $\eps = \lambda-\lambda_1$, then
\begin{gather}
C^2\int\limits_{-\infty}^{+\infty} e^{-\lambda |x|}e^{-\lambda_1 |x-t|}dx 
\le C^2\int\limits_{-\infty}^{+\infty} e^{-\lambda |x|}
e^{-\lambda_1 (|t|-|x|)}dx \notag \\
\le C^2\int\limits_{-\infty}^{+\infty} e^{-\eps |x|} e^{-\lambda_1 |t|}dx \le
\le C^2 e^{-\lambda_1 |t|} \int\limits_{-\infty}^{+\infty} e^{-\eps |x|}dx 
\notag
\\ \Rightarrow |(\stft _g g)(t,\omega)| \le C_1  e^{-\lambda_1 |t|}
\label{stft1}
\end{gather}
for some constant $C_1$ depending on $\lambda_1$.

Similarly 
\begin{equation}
\label{stft2}
|(\stft _g g)(t,\omega)| = |(\stft_{\hat{g}} \hat{g})(\omega,-t)| \le 
D_1 e^{-\lambda_1 |\omega|}.
\end{equation}
By combining~\eqref{stft1} and~\eqref{stft2} and taking the square root
we get
\begin{equation}
\label{stft3}
|(\stft_g g)(t,\omega)|\le\sqrt{C_1 D_1}e^{-\frac{\lambda_1}{2}(|t|+|\omega|)}.
\end{equation}
This together with equations~\eqref{innprod} and~\eqref{sampledstft}
yields \eqref{matexp}.

Now consider
\begin{gather}
\|\gamma-\gamman\| \le \|\Ha \HHI \sigma - \Ha \HHI \HHn \HHnI \sigman\| + 
\notag \\ 
\|\Ha \HHI \HHn \HHnI \sigman - \Ha \HHI H \Hna \HHnI \sigman \| \notag \\
\le \underbrace{\|\Ha \HHI (\sigma-\sigman)\|}_{0} + 
A^{-1} \|(\Hn \Hna - H \Hna)\HHnI \sigman\| \notag \\ =
A^{-1} \|(\Hn \Hna - H \Hna)\HHnI \sigman\|. \notag
\end{gather}
Note that $(\Hn\Hna-H\Hna)$ is a matrix that has a finite number of columns
and a biinfinite number of rows, of which the rows with index $|kl| \le n$
are zero. 

By~\eqref{matexp} the entries of $\HHn$ and the nonzero rows of
$(\Hn\Hna-H\Hna)$ decay exponentially off the diagonal.
Relation~\eqref{rieszbounds} implies
\begin{equation}
\cond(\HHn) \le \cond(\HH) \quad \text{for $n=0, 1, \dots$}.
\label{condest}
\end{equation}
Using \eqref{condest} and Proposition~2 in \cite{Jaf90} it follows that
there exists a $\lambda_2$ and a constant $C_2$ depending on $\lambda_2$ 
and on the condition number of $\HH$, such that
$$|[\HHnI]_{k,l;0,0}| \le C_2 e^{-\lambda_2 (|k|+|l|)},$$
and the same is true for $v^{(n)}:=(\Hn\Hna-H\Hna) [\HHnI]_{k,l;0,0}$.
Now it is easy to see that there exist $\lambda'$ and $C'$ such that
$$\|v^{(n)}\|=\sum_{|k|,|l| > n} |v^{(n)}_{kl}|^2 \le C' e^{-\lambda' n}.$$
\end{proof}

For more results on the connection between the decay of a function and its
short time Fourier transform the reader is referred to~\cite{GZ99}.

\bigskip

To fully appreciate the simplicity and advantages of the proposed
approach we consider for comparison the following method for determining 
the ``dual'' expansion coefficients $\langle f , \gamnm  \rangle$.

The function $f \in \LtR$ can be expressed as
\begin{gather}
\label{recon2}
f = \sum_{k,l \in \Zst} \cnm \gnm \\
\text{where}\,\,\, (T T^{\ast})^{+}c=d\,,
\end{gather}
with $d=\{\dnm\}=\{\langle f, \gnm \rangle\}=Tf$.
We can identify $T T^{\ast}$ with its Gram matrix representation
$(T T^{\ast})_{m,n}= \langle \gnmp, \gnm  \rangle$ for $m,n,m',n' \in \Zst.$

Setting $\Tn \Tna = \Pn T T^{\ast} \Pn$ and $\dn = \Pn d$, we
obtain the $n$-th approximation $\cn$ to $c$ by solving
\begin{equation}
\Tn \Tna \cn = \dn.
\label{finitesec}
\end{equation}

Unfortunately the (generalized) inverse of $\Tn \Tna$ is not
bounded for $\ntoinf$, although $\|(T T^{\ast})^{+}\|\le A^{-1}$. In fact, 
$\|(\Tn \Tna)^{+}\| \toinf$ for $\ntoinf$. Hence $\cn$ does not converge to 
$c$ and without further modifications this approach does not lead to an 
approximation of $f$.

Instead of computing $(\Tn \Tna)^{+}$ we can compute a
regularized inverse via a truncated singular value decomposition
by setting the singular values of $T T^{\ast}$ below a certain threshold 
$\tau$ to zero. Let $(\Tn \Tna)^{\tau,+}$ denote this regularized
inverse. It is shown in \cite{Har98} that $(\Tn \Tna)^{\tau,+}$
converges strongly to $(T T^{\ast})^{+}$ if we allow the threshold
parameter to vary for each $n$. 

In order to use this approach for practical purposes we need good
estimates for the optimal threshold $\tau_n$. Assuming a numerical
precision of $\delta$ of the data and setting $B_n = \|\Tn \Tna\|$, it is 
shown in \cite{Str99b} that $\tau_n$ can be estimated by
$$\tilde{\tau}_n \le 
 B_n \left(\frac{\delta}{p}\right)^{\frac{1}{p+1}}.$$
Here -- without going into details about regularization theory --
$p$ can be seen as ``smoothness parameter''~\cite{EHN96}, the
standard setting for $p$ in regularization theory is $p=1$ or $p=2$.
Thus for large $n$ we get
$$\tilde{\tau} \in  [B (\delta/2)^{\frac{1}{3}}, B\delta^{\frac{1}{2}}]$$
where $B$ is the upper frame bound. Good estimates for the upper frame bound 
are important to apply this method in practice.

All this extra effort is not necessary when using Theorem~\ref{th1}. 
A different approach has been proposed in~\cite{CC98}, 
which however requires substantially more effort compared to
the Wexler-Raz based approach in this section.

\section{Laurent operators and decay of dual Gabor frames} \label{s:laurent}

A natural question for Gabor frames -- also in spite of
Theorem~\ref{th:decay} is the following. Given a window $g$ with certain 
decay properties in time and/or frequency, does its dual $\gamma$ have the 
same decay properties? This question is not only of interest from
a theoretical viewpoint, but has a number of practical implications,
e.g., see Section~\ref{s:ofdm}.

We briefly summarize a few important results. We assume in the sequel that 
$\{\gnm\}$ constitutes a Gabor frame for $\LtR$ or $\ltZ$.\\
(i) If $g \in \compact$, then $\gamma \in \compacti$ only in very special 
cases~\cite{boel97c}. In general $\gamma$ is no longer compactly supported, 
but has exponential decay (see~\cite{Str98a} for a proof in $\ltZ$, the
result can be easily extended to $\LtR$). 
Hence in this case $g$ and $\gamma$ do not belong to the same type of
space.\\
(ii) If $g \in \expmat$ then $\gamma \in \expmati$, but in general with a 
different exponent $\gamma_1 < \gamma$~\cite{Str98a,BJ00}. Thus $g$ and 
$\gamma$ have the same type of decay, but do not belong to the same space, 
we lose some quality of decay.\footnote{Note that there do exist Gabor 
frames with Gaussian decay, that have {\em non-canonical} duals with 
Gaussian decay, see Example~3.10 in~\cite{BJ00}.} \\
(iii) If $g \in \polmat$, then $\gamma \in \polmat$ (see~\cite{Str98a} 
for a proof for $\ltZ$, the result can be easily generalized to $\LtR$ 
using the same approach as in~\cite{BJ00}). In this case $g$ and
$\gamma$ actually belong to the same space $\polmat$.

From an algebraic point of view case (iii) is the most appealing one.
Can we find a (simple) condition on the decay of $g$ that implies that
the dual $\gamma$ belongs to the same space, similar to case~(iii)?
In this section we will give an exhaustive answer to this question.

We need some preparation before we proceed.
\begin{definition}
\label{wiener}
Let $w(t)$ be a continuous weight function on $\Rst$. 
The (weighted) Wiener Algebra $\wienerw$ is the Banach space of absolutely 
convergent Fourier series of period 1 (cf.~\cite{Rei68}), i.e., 
$f \in \wienerw$ if 
\begin{equation}
f(t) = \sum_{k =-\infty}^{\infty} c_k e^{2\pi i kt}
\label{fouriersum}
\end{equation}
with
\begin{equation}
\sum_{k =-\infty}^{\infty} |c_k| w(k) < \infty .
\label{finsum}
\end{equation}
The norm on $\wienerw$ is
$$\|f\|_{\wienerw} = \sum_{k=-\infty}^{\infty}|c_k| w_k .$$
\end{definition}
It follows from Chapter 19.4 of~\cite{GRS64} that $\wienerw$ is a Banach
algebra under pointwise multiplication. $\wienerw$ can be identified with 
the space of all sequences $c=\{c_k\}$ which are in $\lsp_{1,w}$.

Definition~\ref{wiener} can be extended to matrix-valued functions in
a straightforward manner. Put
\begin{equation}
\label{definemat}
\Phi(e^{2\pi it})=
\begin{bmatrix}
f_{11}(e^{2\pi it}) & \dots & f_{1m}(e^{2\pi it}) \\
\vdots              &       & \vdots              \\
f_{m1}(e^{2\pi it}) & \dots & f_{mm}(e^{2\pi it}) \\
\end{bmatrix}
\end{equation}
and set
\begin{equation}
\label{fourierint}
A_k = \int \limits_{0}^{1} \Phi(e^{2\pi it}) e^{-2\pi it} dt .
\end{equation}
Let the weight function $w(t)$ act on $\Rst^m$. Then $\Phi$ belongs to
the matrix Wiener algebra $\wienerwm$ if
\begin{equation}
\label{wienersum}
\sum_{k=-\infty}^{\infty} |A_k| w(k) < \infty .
\end{equation}

$\Phi$ is unitarily equivalent to the block Laurent operator $L$
whose matrix representation (with respect to the standard basis)
is given by
$$
\begin{bmatrix}
\ddots &     &        &        &        \\
       & A_0 & A_{-1} & A_{-2} &        \\
       & A_1 & \fbox{$A_0$}  & A_{-1} &        \\
       & A_2 & A_{1}  & A_{0}  &        \\
       &     &        &        & \ddots \\
\end{bmatrix}.
$$
Here \fbox{$A_0$} denotes the $(0,0)$ entry which acts on the $0$-th
coordinate space. $\Phi$ is also called the defining function of the 
block Laurent operator $L$. By a slight abuse of notation we also
write $L=[A_{kl}]_{k,l=-\infty}^{\infty}$ where $A_{kl}=A_{k-l}$.

We define $\laurentwm$ as the space consisting of all block 
Laurent operators whose $m \times m$ blocks $A_k$ satisfy~\eqref{wienersum}.
If $m=1$ we simply write $\laurentw$ and $L$ reduces to a scalar-valued 
Laurent operator.

\begin{theorem}
\label{th:gelfand}
Let $L \in \laurentw$ be self-adjoint and positive definite.
Assume that the weight function $w$ satisfies 
\begin{equation}
\lim_{\ntoinf} \frac{1}{^n\!\sqrt{w(-n)}} = 1
\quad \text{and}\,\,\,
\lim_{\ntoinf} \,\, ^n\!\sqrt{w(n)} = 1 ,
\label{GRS}
\end{equation}
then $L^{-1} \in \laurentw$.
\end{theorem}

\begin{proof}
Set $L=[A_{kl}]_{k,l=-\infty}^{\infty}$ and 
$L^{-1}=[B_{kl}]_{k,l=-\infty}^{\infty}$. Since $L$ is positive definite we 
have
\begin{equation}
f(\omega) = \sum_{k=-\infty}^{\infty} A_k e^{2\pi i k \omega} > 0
\label{}
\end{equation}
and by the properties of Laurent operators~\cite{GGK93}
$$\frac{1}{f(\omega)} = \sum_{k=-\infty}^{\infty} B_k e^{2\pi i k \omega}.$$
The property $L \in \laurentw$ is equivalent to $f \in \wienerw$.
By Theorem~2 on page~24 in~\cite{GRS64} an element of $\wienerw$ has an
inverse in $\wienerw$ if it is not contained in a maximal ideal.
Any maximal ideal of $\wienerw$ consists of elements of the form
(cf.~Chapter 19.4 in~\cite{GRS64})
\begin{equation}
\sum_{k=-\infty}^{\infty}a_k \xi^k = 0, \notag
\end{equation}
where $\xi = \rho e^{2\pi i \omega}$ with
\begin{equation}
\rho_1 \le \rho \le \rho_2 \notag
\end{equation}
and
\begin{equation}
\rho_1 = \underset{\ntoinf}\lim \frac{1}{^n\!\sqrt{w_{-n}}}
\,\,\,\,\text{and} \,\,\,
\rho_2 = \underset{\ntoinf}\lim \, ^n\!\sqrt{w_{n}}.
\notag
\end{equation}
Due to assumption~\eqref{GRS} we get $\rho_1=\rho_2 =1$,
hence $\rho=1$. Thus a necessary and sufficient condition for an
element in $\wienerw$ to be not contained in a maximal ideal 
is $\sum_{k} a_k e^{2\pi i k \omega} \neq 0$ for all $\omega$.
By assumption $L$ is positive definite, hence
$f(\omega)=\sum_{k}A_k e^{2\pi i k \omega}>0$ for all $\omega$,
consequently $L^{-1} \in \laurentw$.
\end{proof}

The proof of Theorem~\ref{th:gelfand} is essentially based on results of
Gelfand, Raikov, and Shilov~\cite{GRS64}. A crucial role plays 
condition~\eqref{GRS}, to which we will henceforth refer as
the {\em Gelfand-Raikov-Shilov condition} (GRS-condition for short) .

\medskip

Now let $L$ be a block Laurent operator with defining function
$\Phi=[f_{ij}]_{i,j=1}^m \in \wienerwm$. If $L$ is hermitian, positive 
definite then
$$\det ([f_{ij}(\lambda)]_{i,j=1}^m) > 0, \qquad |\lambda|=1.$$
Hence we can apply Theorem~8.1 on page~830 in~\cite{GGK93} to extend
Theorem~\ref{th:gelfand} to block Laurent operators and obtain the following
\begin{corollary}
\label{cor:laurent}
Let $L \in \laurentwm$ be self-adjoint and positive definite.
Assume that the weight function $w$ defined on $\Rst^m$ satisfies 
the GRS-condition, i.e.,
\begin{equation}
\lim_{\ntoinf} \frac{1}{^n\!\sqrt{w(-n)}} = 1
\quad \text{and}\,\,\,
\lim_{\ntoinf} \, ^n\!\sqrt{w(n)} = 1 ,
\label{GRSblock}
\end{equation}
then $L^{-1} \in \laurentwm$.
\end{corollary}

It is clear that in the derivations above we can replace 1-periodic
functions by $\alpha$-periodic functions ($\alpha \in \Qst$).
Similarly, sequences with indices in $\Zst$ can be replaced by sequences
indexed by $\Zst_{\alpha}$ (i.e. by indices of the form $k \alpha$).

Now we are ready to prove a general result on the decay properties of
Gabor windows $g$ and their duals $\gamma$.

\begin{theorem}
\label{th:decay}
Let $(g,a,b)$ generate a Gabor frame for $\LtR$ with $ab=\frac{p}{q}\in
\Qst$ with relative prime integers $p$ and $q$. Let $g \in \Lsp_{1,w}(\Rst)$ 
where $w$ satisfies the GRS-condition. 
Then $\gamma \in \Lsp_{1,w}(\Rst)$.
\end{theorem}
\begin{proof}
The assumption $g \in \Lsp_{1,w}(\Rst)$ implies $\sum_{k=-\infty}^{\infty}
|g(k)| w(k) < \infty$ and $\sum_{k=-\infty}^{\infty}|g(k/b)|w(k/b) < \infty$
for $b \in \Nst$. Denote
\begin{equation}
\label{corr}
G_{kl}(t)=\frac{1}{b}\sum_{r=-\infty}^{\infty}g(t-ra-\frac{k}{b})
g^{\ast}(t-ra-\frac{l}{b}), \qquad k,l \in \Zst.
\end{equation}
$G(t)$ is periodic in $t$ with period $a$. Since $ab=\frac{p}{q}$
we also get $G_{kl}(t) = G_{k+p,l+p}(t)$, cf.~\cite{FG97}, Theorem~3.2. 
In words, $G(t)$ is a block Laurent operator. The frame property implies 
that $G(t)$ is hermitian positive definite and
$$A I \le G(t) \le BI.$$

The submultiplicativity of the function $w$ implies that the spaces
$\Lsp_{1,w}(\Rst)$, $\lsp_{1,w}(\Zst)$, and $\lsp_{1,w}(\Zst_{\frac{1}{b}})$
are Banach algebras under convolution (cf.~\cite{GRS64}).
Hence $\{G_{0k}(t)\}_{k \in \Zst} \in \lsp_{1,w}(\Zst^m_{\frac{1}{b}})$,
and $G(t) \in \laurentwm$. By Corollary~\ref{cor:laurent} 
\begin{equation}
\label{inv}
G^{-1}(t)\in \laurentwm.
\end{equation}
The dual $\gamma$ can be expressed as (see e.g.~\cite{BJ00})
\begin{equation}
\gamma(t) = \sum_{k=-\infty}^{\infty}[G^{-1}(t)]_{0k} g(t-\frac{k}{b}).
\label{dual}
\end{equation}
This equation together with~\eqref{inv} and the assumption 
$g \in \Lsp_{1,w}(\Rst)$ yields that $\gamma \in \Lsp_{1,w}(\Rst)$.
\end{proof}

\remark
(i) The idea to exploit the Laurent operator property of the Gabor frame
operator in the context of window decay is not new. It has been used by 
Feichtinger and Gr\"ochenig~\cite{FG97} in connection with
modulation spaces and polynomial weights.\\
(ii) It is easy to reformulate Theorem~\ref{th:decay} for windows whose 
decay properties are given in the frequency domain. We leave this
modification to the reader.

Observe that condition~\eqref{GRS} is satisfied e.g.\ for
$w(x)=(1+|x|)^s$ and $w(x)=e^{\lambda |x|^\gamma}, \gamma<1$, but
not for $w(x)=e^{\lambda |x|}$. This is why we have to use
a smaller exponent $\lambda_1 < \lambda$ for the dual to bound the decay
in case of an exponentially decaying window. 

A detailed discussion on the connection between decay properties of 
windows and modulation spaces can be found in~\cite{Gro00}.

\section{Orthogonal frequency division multiplexing and Gabor systems}
\label{s:ofdm}
Orthogonal frequency division multiplexing (OFDM) has attracted a great
deal of attention as an efficient technology for wireless data
transmission~\cite{Zou95}. Among others it is currently used in the 
European digital audio broadcasting standard~\cite{CFH89}.
In a wireline environment OFDM is known under the name
{\em discrete multitone transmultiplexing} (DMT).

\begin{figure}
\begin{center}
\epsfig{file=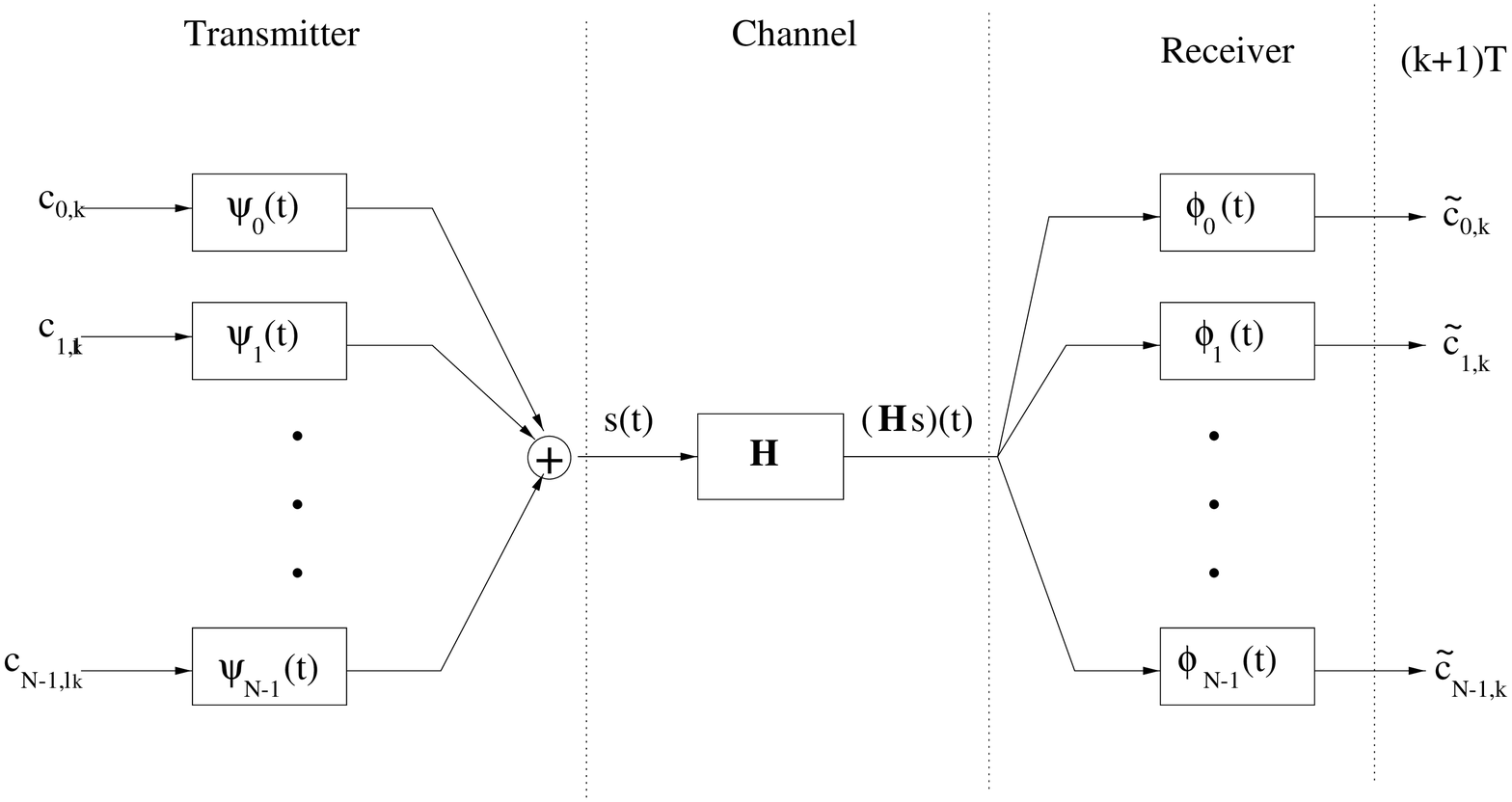,width=130mm}
\caption{A general baseband OFDM system configuration}
\label{fig:ofdm}
\end{center}
\end{figure}

The basic idea of OFDM is to divide the available spectrum into several
subchannels (subcarriers). A baseband OFDM system is schematically represented 
in Figure~\ref{fig:ofdm}.\\
{\bf Transmitter:} Assuming $N$ subcarriers, a bandwidth of $W$~Hz, 
{\em symbol length} of $T$ seconds, and {\em carrier separation} $F:=W/N$, 
the transmitter of a general OFDM system uses the following waveforms 
\begin{equation}
\label{ofdmsystem}
\psi_l(t)=\psi(t) e^{2\pi it lF}, \quad k \in \Zst, l = 0,\dots,N-1.
\end{equation}
The transmitted baseband signal for OFDM symbol number $k$ is
\begin{equation}
\label{OFDMsingle}
s_k(t)=\sum_{l=0}^{N-1} c_{kl} \psi_l(t-kT),
\end{equation}
where $c_{k0}, c_{k1},\dots,c_{k,N-1}$ are the complex-valued information
bearing coefficients (data symbols). Assuming an infinite sequence of
OFDM symbols is transmitted, the output from the transmitter is a
superposition of individual OFDM symbols
\begin{equation}
\label{OFDMsignal}
s(t)=\sum_{k=-\infty}^{\infty} s_k(t) =
\sum_{k=-\infty}^{\infty}\sum_{l=0}^{N-1} c_{k,l} \psi_l(t-kT).
\end{equation}
\noindent
{\bf Receiver:}
The OFDM receiver consists of a matched filter bank $\{\phi_l\}$ of similar
structure as the transmitter waveforms, i.e.,
\begin{equation}
\label{ofdmrec}
\phi_{kl}(t)=\phi(t-kT) e^{2\pi itlF}, \quad k \in \Zst, l = 0,\dots,N-1.
\end{equation}
The transmitted data are recovered by projecting the received signal
$r={\mathbf H}s+\nu$ (where $\nu$ represents additive white Gaussian noise,
AWGN for short) onto the functions $\phi_{kl}$, i.e., 
\begin{equation}
\label{recon}
\tilde{c}_{kl}=\langle r ,\phi_{kl} \rangle.
\end{equation}
In the standard OFDM setup the functions $\psi_{kl}$ are designed
to be mutually orthogonal. In this case $\phi\equiv \psi$. The situation 
where the sets $\{\psi_{kl}\}$ and $\{\phi_{kl}\}$ are biorthogonal
is referred to as {\em biorthogonal frequency division multiplexing} (BFDM).

It is useful to have the following two OFDM setups at hand.

\noindent
{\bf Continuous-time model:}
Since $s(t)$ is a continuous signal, it is convenient to consider a 
continuous-time transmission set that is infinite in both indices $k$ and $l$, 
which means we set $c_{kl}=0$ if $l<0$ or $l>N-1$. Using the viewpoint it
is clear that OFDM can be interpreted as an orthonormal Gabor system in
$\LtR$. For a theoretical analysis it is often
advantageous to work with this time-continuous model.

\noindent
{\bf Discrete-time model:}
In practice the OFDM system is digitally implemented. At the transmitter
the signal $s$ is passed through a digital-to-analog converter to
obtain a continuous signal.
At the receiver an analog-to-digital converter transforms the received
signal back into a discrete-time signal.
Therefore it is useful to consider the following discrete-time setup for OFDM.

Let $\psi \in \ltZ$ and set $\psikl(n)=\psi(n-kT)e^{2\pi ilF}$ for
$k,n \in \Zst, l =0,\dots, N-1$ and $F=W/N$. Then the OFDM signal becomes
\begin{equation}
\label{disctransmit}
s(n)=\sum_{k=-\infty}^{\infty}\sum_{l=0}^{N-1} c_{k,l} \psikl .
\end{equation}
It is clear that $\{\psikl\}$ coincides with a discrete-time
Gabor system.

\subsection{Wireless channels and time-frequency localization} 
\label{ss:channel}

In the ideal case the channel {\bf H} does not introduce any distortion, and
the data $c_{kl}$ can be recovered exactly. However in practice
wireless channels introduce time dispersion as well as frequency
dispersion (in addition to the usual channel noise).
The time dispersion is caused by multipath propagation and can lead to
intersymbol interference (ISI). Frequency dispersion of the mobile radio
channel is due to the Doppler effect and can cause interchannel
interference (ICI).
The distortion resulting from channel dispersion depends crucially
on the time-frequency localization of the transmitter pulse shapes
$\{\phi_{kl}\}$. Robustness against doubly dispersive channels can be
achieved by pulse shapes with good time-frequency localization.
 
An optimum OFDM system in case of doubly dispersive channels would
consist of orthogonal basis functions with $TF=1$, such that
the $\psikl$ are well localized in time and frequency.
The condition $TF=1$ (critical sampling) ensures maximal spectral
efficiency of the transmission system. Unfortunately, such a system cannot 
exist due to the Balian-Low theorem~\cite{Dau92}.

Therefore other approaches (cyclic prefix, 
pulse-shaping, BFDM) have been proposed for doubly dispersive
channels~\cite{Zou95,HB97,KM98,Bol99}. These approaches can be interpreted
in the time-frequency plane as using an undersampled grid, i.e., $TF>1$,
which results in a set of basis functions that is incomplete
in $\LtR$ (or $\ltZ$, respectively).
Although the choice $TF>1$ leads to a loss in the capacity of the transmission 
system (determined by the undersampling rate), it is usually an acceptable 
price to pay to mitigate interference.

We know from Gabor theory that there exist Gabor systems with good 
time-frequency localization for $TF>1$ (take for instance $TF=2$ and a 
Gaussian window). In the present context of OFDM we are particularly
interested in the construction of orthonormal incomplete Gabor systems
with good time-frequency localization properties.

The following result ties tight Gabor frames for $\LtR$ to orthonormal
Gabor bases for subspaces of $\LtR$. It is a simple consequence of 
Theorem~\ref{ronshen1}.
\begin{corollary} [\cite{RS97a}]
\label{ronshen2}
Assume $\|g\|=1$. $\gba$ is an orthonormal basis for its closed linear span 
if and only if $\gab$ is a tight frame for $\LtR$.
\end{corollary}

Now, let $\{\gnm\}$ be a frame for $\LtR$ with frame operator $S$, then a 
standard way to construct a tight frame is the following. Compute 
\begin{equation}
\label{ortho}
\gtight = S^{-\frac{1}{2}}g
\end{equation}
then the set $\{\gtightnm\}$ is a tight frame for $\LtR$ and by 
Corollary~\ref{ronshen2} (after normalizing $\gtight$) $\{\gtightkl\}$ is an 
orthonormal basis for its closed linear span (equivalently we can apply this 
orthogonalization procedure directly to the set $\{\gkl\}$). 
Our OFDM system is now given by setting $\psi:=\gtight$, $T=1/b$, $F=1/a$,
or equivalently,
\begin{equation}
\label{gabor2ofdm}
\psi_{kl}:= \gtightkl .
\end{equation}

In order to obtain an OFDM system that is well localized in the time-frequency 
plane, it seems natural to start with a function $g$ with good time-frequency 
localization such that $(g,a,b)$ generates a frame for $\LtR$, apply the
orthogonalization procedure~\eqref{ortho} and hope that 
$\gtight=S^{-\frac{1}{2}}g$ and whence $\psi$ inherit these 
localization properties. This approach is not new, it has already been 
considered in~\cite{KM98,Bol99}, however with very different 
conclusions.

It is stated in~\cite{KM98} that ``{\em such an orthogonalization of the 
pulses, however, is not desirable \dots because the good time-frequency 
localization of the pulses is destroyed by such a transformation}''. 
In contrast in~\cite{Bol99} it is claimed that the orthogonalization 
procedure~\eqref{ortho} ``{\em in practice \dots starting from a 
well-localized initial filter \dots yields well-localized 
orthogonal filters}''.

Who is right and who is wrong? The answer is: 
\begin{quote}
{\em Both and nobody -- it depends!}
\end{quote}

The first step to this answer is the following result.

\begin{corollary}
\label{cor:tight}
(i) Let $g,\hat{g} \in \Lsp_{1,w}$ and let $(g,a,b)$ generate a Gabor frame
for $\LtR$. Set $\gtight = S^{-\frac{1}{2}} g$. If the weight function
$w$ satisfies the GRS-condition, then $\gtight,\hat{\gtight}\in \Lsp_{1,w}$.\\
(ii) Let $g, \hat{g} \in \lsp_{1,w}$ and let $(g,a,b)$ generate a Gabor frame
for $\ltZ$. Set $\gamma = S^{-1}g$ and $\gtight = S^{-\frac{1}{2}} g$. 
If the weight function $w$ satisfies the GRS-condition, then 
$\gamma, \hat{\gamma} \in \lsp_{1,w}$ and
$\gtight,\hat{\gtight} \in \lsp_{1,w}$.
\end{corollary}

\begin{proof}
(i) The proof of this result is similar to the proof of
Theorem~\ref{th:decay}. We only indicate the necessary  modifications.
Let $\Phi$ be the matrix-valued defining function of a hermitian
positive definite Laurent operator $L$ and denote 
$L^{-\frac{1}{2}} = [C_k]_{k,l=-\infty}^{\infty}$.  
It follows from the basic properties of block Laurent operators that 
$\Phi(\omega) > 0$ for all $\omega$ and
$$C_k=\int\limits_{0}^{1}\frac{1}{\sqrt{\Phi(\omega)}} 
\,\, e^{-2\pi ik\omega}\, d\omega , \qquad k \in \Zst.$$
The fact that $L$ is self-adjoint positive-definite implies that
$L^{\frac{1}{2}}$ is in the same algebra as $L$, see~\cite{Gar66}.
Thus if $L \in \laurentwm$ (and consequently $L^{\frac{1}{2}}\in \laurentwm$)
for weights satisfying~\eqref{GRSblock}, then
$L^{-\frac{1}{2}} \in \laurentwm$. The rest follows now by repeating the steps
in Theorem~\ref{th:decay} and using the remark following
Theorem~\ref{th:decay}.\\
(ii) In this case the frame operator $S$ is a block Laurent operator
(see~\cite{Str98a}). The assumption $g \in \lsp_{1,w}$ implies that 
$S \in \laurentwm$. Along the same lines as above and in
Theorem~\ref{th:decay} it follows that $S^{-1}$ and 
$S^{-\frac{1}{2}} \in \laurentwm$ and consequently 
$\gamma,\hat{\gamma} \in \lsp_{1,w}$, similar for $\gtight,\hat{\gtight}$.
\end{proof}
Thus if the window (i.e.\ pulse shape) $g$ satisfies the properties of
Corollary~\ref{cor:tight} the orthogonalization procedure~\eqref{ortho}
yields a function that is in the same space (in the same {\em algebra})
as the function $g$, and the quality of decay does not 
change.\footnote{Due to the
connection between Gabor frames and modulated filter banks~\cite{CV98a},
Corollary~\ref{cor:tight} also yields a characterization of the
decay properties of paraunitary modulated filter banks.}

For a window $g$ with exponential decay (in time and frequency) it has been 
shown that $\gtight$ also has exponential decay (see~\cite{Str98a} for a 
proof for $\ltZ$, this result has later been extended to $\LtR$, 
cf.~\cite{BJ00}). However, similar to the dual window, in general the 
exponent for the decay of $\gtight$ will be smaller than for $g$. 

Since the Gaussian is optimally localized in the time-frequency plane,
in the sense that it minimizes the uncertainty principle, it is
interesting to note that for $g(t)=e^{-\pi t^2}$ the function
$\gtight$ has only exponential decay in time and frequency, cf.~\cite{BJ00}.
Thus we certainly loose some time-frequency localization in this case, 
although the resulting function $\gtight$ still has exponential decay. 

Thus it seems there exists a large class of windows whose time-frequency
localization properties are not affected by applying~\eqref{ortho}.
This is however only half of the truth, since in the considerations above 
we have ignored any constants that come into play. 
This is acceptable from an asymptotic-analysis viewpoint,
but not for applications.
 
Consider for example the following situation. Take a window $g$ with
exponential decay in time and frequency and assume that
$(g,T,F)$ constitutes a frame for $TF>1$.
Then for any $\eps >0$ there exists an $N$ such that for $T>N, F>N$
$$\Big\|\frac{g}{\|g\|}-\frac{\gtight}{\|\gtight\|}\Big\| \le \eps$$ 
since $(TF) \gamma \rightarrow g$ for $T\rightarrow\infty,F\rightarrow \infty$
(see~\cite{FZ98} for a mathematically precise formulation).
But of course $TF \gg 1$ leads to an unacceptable large loss of capacity for 
OFDM, already $TF>2$ seems to be prohibitive in this context.

On the other hand, if we let $TF\rightarrow 1$, then $\gtight$ and
$\gamma$ will become increasingly ``ill-localized'' in the time-frequency
plane (since in the limit case $TF=1$ we are confronted with the Balian-Low
theorem), although from theory we know that $\gtight$ has 
exponential decay as long as $TF>1$. However, as pointed out earlier, this 
decay involves a constant $C$ that depends on the ratio of the frame bounds.
Since $A\rightarrow 0$ for $TF\rightarrow 1$ it follows that 
$C \rightarrow \infty$. Thus for applications such as OFDM the 
exponential-decay property of $\gamma$ and $\gtight$ quickly becomes 
meaningless for $TF$ close to 1. If we choose a well-localized window,
there is a trade-off between increasing $TF$ and increasing the frame
bound ratio $B/A$. 

Hence a correct formulation from a practical viewpoint of the two contradicting 
statements above is:
\begin{quote}
{\em If $g$ is well localized in time and frequency and if the frame
bounds satisfy $B/A \approx 1$, then $\gtight$ is also well localized in
time and frequency.}
\end{quote}

Numerical experiments indicate that for $TF=1.3$ or $TF=1.4$ it is possible
to construct OFDM basis functions with good time-frequency localization,
e.g.\ choose the Gaussian as initial window $g$. A definite answer if the
time-frequency properties of the resulting OFDM basis functions obtained in
that way are sufficient for practical purposes is difficult, since it depends
on the actual ISI and ICI, as well as the AWGN behavior of the channel.

\bigskip
\remark 
It is well-known that Wilson bases can be constructed from twofold
oversampled tight Gabor frames~\cite{DJJ91}. In light of this
fact Corollary~\ref{cor:tight} provides an extension of existing results on
the connection between the decay behavior of Gabor frames and
Wilson bases.

\end{document}